\documentclass[a4paper,10pt,draft]{article}

\usepackage[latin1]{inputenc}
\usepackage{amsfonts}
\usepackage{amsmath}
\usepackage{amssymb}
\usepackage{amsthm}
\usepackage[french]{babel}
\usepackage[T1]{fontenc}
\usepackage[dvips]{graphicx}
\usepackage[all]{xy}

\newtheorem*{definition}{Définition}
\newtheorem*{proposition}{Proposition}

\theoremstyle{definition}

\theoremstyle{remark}
\newtheorem*{remarque}{Remarque}

\newcommand{\opun}{_{(1)}}
\newcommand{\opn}{_{(n)}}
\newcommand{\opnp}{_{(n+1)}}
\newcommand{\opz}{_{(0)}}
\newcommand{\opmun}{_{(-1)}}
\newcommand{\optr}{_{(3)}}
\newcommand{\opmd}{_{(-2)}}
\newcommand{\opd}{_{(2)}}
\newcommand{\opzt}{_{\tilde{(0)}}}
\newcommand{\opmt}{_{(-3)}}
\newcommand{\opq}{_{(4)}} 

\newcommand{\BZ}{\Bbb{Z}}

\newcommand{\lra}{\longrightarrow}

\begin{document}

\centerline{\bf ALG\`EBRES DE FROBENIUS - VIRASORO}

\vspace{0.25in}

\centerline{Maxime \textsc{Rebout}\footnote{Laboratoire Emile Picard, 
Universit\'e Paul Sabatier, 118 Route de Narbonne, 31062 Toulouse}, 
Vadim \textsc{Schechtman}\footnotemark[1]} 

\vspace{1in}

\section{Introduction}

Consid\'erons un couple $\frak{g} = (\frak{l}, \langle, \rangle)$, o\`u 
$\frak{l}$ est une alg\`ebre de Lie sur un anneau commutatif de base $k$ 
fix\'e et 
$$
\langle,\rangle:\ \frak{l}\times \frak{l} \lra k 
$$
est une forme $k$-bilin\'eaire sym\'etrique $\frak{l}$-invariante 
(on pourrait appeler un tel couple {\em une alg\`ebre de Lie Frobenius}). 
On supposera d\'esormais que $k$ contient $\Bbb{Q}$.  
On sait associer \`a $\frak{g}$ une alg\`ebre vertex $V_{\frak{g}}$, 
le module "vacuum" sur l'alg\`ebre de Kac-Moody $\hat{\frak{g}}$ correspondante,
[K]. 

Celle-ci est une $k$-alg\`ebre vertex $\BZ_{\geq 0}$-gradu\'ee, 
$V_{\frak{g}} = \oplus_{n\geq 0}\ V_{\frak{g},n}$ avec $V_{\frak{g},0} = 
k\cdot 1$, $1$ d\'esignant le vecteur vacuum, $V_{\frak{g},1} = \frak{l}$, 
la structure d'un $k$-module sur $V_{\frak{g},1}$ est induite par l'op\'eration 
$\opmun:\ V_{\frak{g},0}\times V_{\frak{g},1} \lra V_{\frak{g},1}$,
l'op\'eration 
$$
\opz:\ V_{\frak{g},1}\times V_{\frak{g},1} \lra V_{\frak{g},1}
$$
co\"\i ncide avec le crochet de Lie sur $\frak{l}$ et l'op\'eration
$$
\opun:\ V_{\frak{g},1}\times V_{\frak{g},1} \lra V_{\frak{g},0}
$$
co\"\i ncide avec l'accouplement $\langle, \rangle$. 

Alors $V_{\frak{g}}$ est caract\'eris\'ee par la propri\'et\'e universelle 
suivante: un morphisme de $V_{\frak{g}}$ dans une alg\`ebre vertex  
$V$ quelconque est la m\^eme chose qu'un morphisme de 
$k$-modules $\phi:\ \frak{l} \lra V$ tel que 
$$
\phi([x,y]) = \phi(x)\opz\phi(y)\ 
$$
et
$$ 
\langle x,y\rangle \cdot 1_V = \phi(x)\opun\phi(y),\ 
x, y\in \frak{l}.
$$
En effet, une alg\`ebre de Lie Frobenius $\frak{g}$ est un exemple simple
d'une {\em alg\'ebro\"\i de vertex}, 
et l'alg\`ebre vertex $V_{\frak{g}}$ est son alg\`ebre enveloppante, [GMS]. 

Maintenant partons d'une {\em alg\`ebre commutative de Frobenius} 
(par la suite, on écrira plus simplement alg\`ebre de Frobenius). Rappelons 
que c'est une $k$-alg\`ebre commutative associative unitaire $F$, d'unit\'e 
$e$, munie d'une forme $k$-bilin\'eaire sym\'etrique 
$$
\langle , \rangle:\ F\times F \lra k 
$$
telle que, pour tout $x$, $y$, $z$ appartenant à $F$, on ait 
$$
\langle xy, z\rangle = \langle x, yz\rangle,
$$
Nous appelerons le nombre 
$$
c := \langle e, e\rangle
$$
{\em la charge} de $F$. 

{\bf Exemple.} \'Etant donn\'e $c\in k$, l'alg\`ebre de Frobenius $k_c$ co\"\i
ncide avec $k$ comme une alg\`ebre commutative, et 
$\langle x, y\rangle = xyc$. 

Dans cette note on associe, de fa\c{c}on fonctorielle, \`a une 
alg\`ebre commutative de Frobenius $F$ 
une alg\`ebre vertex $\BZ_{\geq 0}$-gradu\'ee
$V_F$. Pour $F= k_c$, $V_F$ est l'alg\`ebre vertex de Virasoro de charge
centrale $c$ (le module "vacuum" sur l'alg\`ebre de Lie de Virasoro), [K]. 
Plus g\'en\'eralement, pour $F$ arbitraire, l'unit\'e $e\in F$ donne lieu 
\`a un vecteur de Virasoro dans $V_F$.  
On pourrait appeler les alg\`ebres 
$V_F$ {\em alg\`ebres de Frobenius - Virasoro}. 

En effet, on d\'efinit, par une m\'ethode similaire \`a celle de [GMS], 
une notion d'une {\em alg\`ebro\"\i de de Virasoro}, dont l'alg\`ebre de
Frobenius est un cas tr\`es particulier. Nos alg\`ebres $V_F$ seront  
les alg\`ebres enveloppantes de ces alg\`ebro\"\i des.

\section{Alg\`ebres vertex du type Virasoro}

On commence par regarder l'exemple de l'algèbre 
vertex de Virasoro. Cette algèbre vertex $\BZ_{\geq 0}$-gradu\'ee 
$V = Vir_c$ (tous les objets considérés vivront au-dessus de  
$k$) est par 
définition engendrée par un vecteur $L$ appartenant à $V_2$ qui vérifie l'OPE 
suivante:
$$
L(z)L(w)\ \sim\ \frac{\frac{1}{2}c}{(z - w )^4} +
\frac{2L(w)}{(z - w)^2} +
\frac{\partial L(w)}{(z - w)}
$$
(un nombre $c\in k$ \'etant fix\'e). 
Ceci est équivalent à:
$$
L\optr L=\frac{c}{2};
$$ 
$$
L\opd L=0;
$$ 
$$
L\opun L=2L;
$$ 
$$
L\opz L=\partial L.
$$
A partir de là, on en déduit que 
$$
V_0=k.1;
$$ 
$$
V_1=0;
$$ 
$$
V_2=k.L;
$$ 
$$
V_3=k.\partial L.
$$

On veut maintenant généraliser cet exemple en considérant les algèbres 
vertex qui possèdent une structure analogue à celle de l'algèbre 
de Virasoro. Soit $A$ un $k$-module; on regarde donc plus 
en détails les algèbres vertex du type suivant:
$$
\xymatrix{ A=V_0 \ar[d]^\partial\\ 
0 \ar[d]^\partial\\ V_2 \ar[d]^\partial\\ 
V_3\simeq\partial V_2 \ar[d]^\partial\\
  \vdots}
$$
Donc $V_1 = 0$ et $\partial$ induit un isomorphisme de $V_2$ sur $V_3$. 
Une telle alg\`ebre vertex sera appel\'ee une alg\`ebre vertex {\em du type 
Virasoro}.  

Soit $V$ une telle algèbre vertex. Par la suite, et sauf indications 
contraires, on désignera les éléments de $A$ par les lettres $a$, $b$, $c$ 
et les éléments de $V_2$ par les lettres $x$, $y$ et $z$. De plus, on abrégera 
la notation de l'opération $\opmun$ en écrivant juste $xy$ à la place de 
$x\opmun y$, pour tout $x$, $y$ dans $V$.
  
\section{Premières remarques}
\label{remar}

L'opération $\opmun$, lorsqu'elle est restreinte à $A$, devient associative et commutative car $a\opn b=0$ pour tout $n\geq 0$; ce qui implique que le module $A$ est muni d'une structure de $k$-algèbre unitaire (d'unité 1).

Par ailleurs, le fait que $V_1=0$ permet d'affirmer la trivialité de toutes les opérations de la forme $a\opn y$ pour tout $y\in V$ et tout $n\neq -1$. En effet, on a:
$$a\opn y=-\frac{1}{n+1}(\partial a)\opnp y=0 \text{ car }\partial a\in V_1.$$

Si on s'intéresse maintenant à l'action de $A$ sur $V_2$, c'est-à-dire à l'opération $V_0\times V_2 \overset{_{(-1)}}{\longrightarrow}V_2$, on s'aperçoit que:
\begin{eqnarray*}
(ab)x & = & abx + b\opmd a\opz x + a\opmd b\opz x + b\opmt a\opun x + a\opmt b\opun x +0\\
& = & abx\qquad\text{d'après la remarque précédente,}
\end{eqnarray*}
et \begin{eqnarray*}
xa &=& ax - \partial(a\opz x) + \frac{\partial^2}{2}(a\opun x)+0\\
&=& ax\qquad\text{toujours d'après la même remarque}.
\end{eqnarray*}
Enfin, on sait que $1x=x$ pour tout $x\in V_2$. Ces trois équations nous permettent donc d'affirmer que l'opération $\opmun$ induit une structure de $A$-module bilatère sur $V_2$.
Par ailleurs, le même raisonnement s'applique si on remplace $V_2$ par $V_3$; le $k$-module $V_3$ possède donc lui aussi une structure de $A$-module bilatère.

Pour terminer ces remarques, il faut voir que l'affirmation $V_3\simeq\partial V_2$ n'est pas absurde: en effet, {\itshape a priori}, l'application $\partial$ n'est pas $A$-linéaire, mais seulement $k$-linéaire; donc il serait plus naturel de penser qu'on peut seulement avoir $V_3\simeq A\partial V_2$. Pour autant, comme $V_1=0$, on a $\partial A=0$ et comme l'application $\partial$ est une dérivation, on peut écrire:
$$\partial(ax)=\partial(a)x+a\partial(x)=a\partial x.$$
Donc, dans ce cas particulier, $\partial$ est $A$-linéaire. On peut 
donc bien supposer que $V_3\simeq\partial V_2$.

\section{Structure de $V_2$}

On veut maintenant regarder plus précisément les opérations présentes 
sur $V_2$ et quelles relations ces opérations vérifient. Les opérations 
qui nous intéressent vont donc être les suivantes: 

\begin{itemize}
\item $V_2\times V_2 \overset{_{(1)}}{\longrightarrow}V_2$;
\item $V_2\times V_2
  \overset{_{(0)}}{\longrightarrow}V_3
  \overset{\partial^{-1}}{\longrightarrow}V_2$, 
cette seconde opération sera notée $\opzt$;
\item $V_2\times V_2 \overset{_{(3)}}{\longrightarrow}V_0$, cette dernière 
opération sera notée $\langle.,.\rangle$.
\end{itemize}

\begin{remarque}
  {\itshape A priori}, il y a une quatrième opération à considérer: l'action de $V_2$ sur $V_0$. Mais en regardant de plus près cette application, on s'aperçoit qu'elle est triviale. En effet: $$x\opun a= a\opun x=0,$$ d'après une des remarques initiales.

\end{remarque}

\subsection{$(1):\ V_2\times V_2\longrightarrow V_2$} 

\begin{eqnarray*}
x\opun y &=& y\opun x-\partial\underbrace{(y\opd x)}_{\in V_1}+ 
\frac{\partial^2}{2}\underbrace{(y\optr x)}_{\in V_0}\\
&=& y\opun x,
\end{eqnarray*}

\begin{eqnarray*}
(x\opun y)\opun z &=& x\opun y\opun z + y\opd x\opz z - x\opz y\opd z -y\opun
x\opun z\\
&=& x\opun y\opun z -y\opun x\opun z+ y\opd x\opz z
\end{eqnarray*}
or $y\opd x\opz z=-(x\opz z)\opd y=-\partial(\partial^{-1}(x\opz z))\opd y=2(x\opzt z)\opun y$, donc $$
(x\opun y)\opun z =x\opun y\opun z -y\opun x\opun z+ 2(x\opzt z)\opun y.$$
L'opération $\opun$ est donc commutative, mais non associative.
$$(ax)\opun y=a(x\opun y)+x\opz a\opz y+a\opmd x\opd y+x\opmun a\opun y + a\opmt x\optr y=a(x\opun y).$$
L'opération $\opun$ possède donc en plus la propriété d'être $A$-bilinéaire.

\subsection{${\tilde{(0)}}:\ V_2\times V_2\longrightarrow V_3\longrightarrow V_2$}
\begin{eqnarray*}
x\opzt y &=& \partial^{-1}(x\opz y)\\
&=& \partial^{-1}(-y\opz x+\partial(y\opun x))\\
&=& -y\opzt x + y\opun x.
\end{eqnarray*}
Cette opération n'est donc pas commutative, mais on déduit de l'équation précédente que la partie symétrique de $\opzt$ est $\frac{1}{2}\opun$.

$A$-bilinéarité:
\begin{eqnarray*}
  (ax)\opzt y & = & \partial^{-1}\left((ax)\opz y\right)\\
  & = & \partial^{-1}\left(a(x\opz y)+x(a\opz y)+a\opmd(x\opun y)+x\opmd(a\opun z)+a\opmt(x\opd y) \right)\\
  & = & \partial^{-1}\left(a(x\opz y)\right)\qquad\qquad\text{d'après la section \ref{remar}}\\
  & = & a\partial^{-1}\left(a(x\opz y)\right)\qquad\text{car l'application }\partial^{-1}\text{ est $A$-linéaire}\\
  & = & ax\opzt y.
\end{eqnarray*}
L'opération $\opzt$ est donc $A$-linéaire par rapport à la première variable. Pour prouver la linéarité par rapport à la seconde variable, on utilise l'équation précédente et la $A$-bilinéarité de l'opération $\opun$. Ici encore, on a bien la $A$-bilinéarité de l'opération étudiée.

Relation d'"associativité":
\begin{eqnarray*}
(x\opzt y)\opzt z &=&\partial^{-1}(\partial^{-1}(x\opz y)\opz z)\\
&=& \partial^{-1}(-(x\opz y)\opun z)\\
&=& -\partial^{-1}(x\opz y \opun z - y\opun x\opz z)\\
&=& -x\opzt y\opun z + \partial^{-1}(y\opun x\opz z)\\
&=& -x\opzt y\opun z + \partial^{-1}((x\opz z)\opun z - \partial((x\opz z)\opd y))\\
&=& -x\opzt y\opun z + \partial^{-1}(-(x\opzt z)\opz y)-(x\opz z)\opd y\\
&=& -x\opzt y\opun z -(x\opzt z)\opzt y + 2(x\opzt z)\opun y.
\end{eqnarray*}

\subsection{$\langle.,.\rangle:\ V_2\times V_2\longrightarrow V_0$}
$$\langle x,y\rangle=x\optr y=y\optr x-\partial(y\opq x)+\dots=y\optr x=\langle y,x\rangle,$$
$$\langle ax,y\rangle=(a\opmun x)\optr y=a\opmun x\optr y+x\opd a\opz y+x\opun a\opun y=a\opmun x\optr y=a\langle x,y\rangle.$$
Ce crochet est donc bilinéaire symétrique. 

\section{Relations entre les opérations}
On recherche maintenant les relations qu'il peut y avoir entre les différentes opérations sur $V_2$.
\begin{eqnarray*}
(x\opzt y)\opun z &=& \frac{1}{2}(x\opun y)\opd z\\
&=& -\frac{1}{2}(x\opz y\opd z-y\opd x\opz z)\\
&=& \frac{1}{2}y\opd x\opz z\\
&=& -\frac{1}{2}(x\opz z)\opd z\\
&=& (x\opzt z)\opun y.
\end{eqnarray*}

Maintenant, si on regarde l'équation d'"associativité" de l'opération $\opun$, en prenant les différentes permutations circulaires en $x$, $y$ et $z$ et en sommant les trois équations obtenues, on trouve:
$$2 Cycle_{x,y,z}(x\opun y\opzt z)=Cycle_{x,y,z}(x\opun y\opun z).$$

Par ailleurs,
\begin{eqnarray*}
\langle x\opzt y,z\rangle &=& (x\opzt y)\optr z\\
&=& -\frac{1}{4}(x\opz y)\opq z\\
&=& -\frac{1}{4}(x\opz y\opq z-y\opq x\opz z)\\
&=& -\frac{1}{4}(x\opz z)\opq y\\
&=& \langle x\opzt z,y\rangle,
\end{eqnarray*}

\begin{eqnarray*}
\langle x\opun y,z\rangle &=& (x\opun y)\optr z\\
 &=& x\opun y\optr z+y\opq x\opz z - x\opz y\opq z - y\optr x\opun z\\
 &=& y\opq x\opz z - \langle y,x\opun z\rangle\\
 &=& -(x\opz z)\opq y - \langle y,x\opun z\rangle\\
 &=& 4\langle x\opzt z, y\rangle-\langle y, x\opun z\rangle.
\end{eqnarray*}

Maintenant, si on prend les différentes permutations circulaires en $x$, $y$ et $z$ de l'équation précédente et en les sommant, on trouve:
$$
Cycle_{x,y,z} \langle x\opun y,z\rangle=2Cycle_{x,y,z} 
\langle x\opzt y,z\rangle.
$$

\section{Relation avec les algèbres de Frobenius} 

Dans les sections précédentes, nous avons donc mis en avant une 
structure particulière sur $V_2$. Appelons {\em une alg\`ebro\"\i de de
Virasoro} un couple $(A = V_0, V_2)$, o\`u $A$ est une $k$-alg\`ebre 
commutative, $V_2$ est un $A$-module muni des trois op\'erations suivantes:
\begin{itemize}
\item $\opun:V_2\times V_2\lra V_2$ qui doit être symétrique $A$-bilinéaire;
\item $\opzt:V_2\times V_2\lra V_2$ qui doit juste être $A$-bilinéaire;
\item $\langle .,.\rangle:V_2\times V_2\lra V_0$ qui doit être $A$-bilinéaire symétrique.\\
\end{itemize}
Par ailleurs, ces opérations doivent vérifier les six égalités suivantes:
\begin{gather}
  (x\opun y)\opun z =x\opun y\opun z -y\opun x\opun z+ 2(x\opzt z)\opun y\\
  x\opzt y = -y\opzt x + y\opun x\\
  (x\opzt y)\opzt z=-x\opzt y\opun z -(x\opzt z)\opzt y + 2(x\opzt z)\opun y\\
  (x\opzt y)\opun z = (x\opzt z)\opun y\\
  \langle x\opzt y,z\rangle = \langle x\opzt z,y\rangle\\
  \langle x\opun y,z\rangle = 4\langle x\opzt z, y\rangle-\langle y, x\opun z\rangle
\end{gather}
On dira aussi que $V_2$ est une alg\'ebro\"\i de de Virasoro 
sur $A$. \\

Les consid\'erations des sections pr\'ec\'edentes fournissent le foncteur 
$$
t:\ \text{(Alg\`ebres\ du\ type\ Virasoro)} \lra 
\text{(Alg\'ebro\"\i des\ de\ Virasoro)},\ 
$$
$$
t(V) = (V_0,V_2)
$$
Ce foncteur admet un adjoint \`a gauche $U$, dit {\em l'alg\`ebre
enveloppante}. 

Si ${\cal{A}}= (A,V_2)$ est une alg\'ebro\"\i de de Virasoro, alors on
construit $U{\cal A}$ en deux \'etapes. D'abord, on d\'efinit \`a partir de ${\cal A}$ 
une alg\`ebre vertex de Lie\footnote{la terminologie de Frenkel - Ben-Zvi; 
dans [K] un tel objet est appel\'e une alg\`ebre conforme}
 $\Bbb Z_{\geq 0}$-gradu\'ee (cf. [GMS], D\'efinition 10.1) 
$$
L{\cal A} = 
A \oplus V_2 \oplus \partial V_2 \oplus \partial^2 V_2 \oplus \ldots 
$$
Ici $\partial^i V_2$ est une copie de $V_2$; on introduit les op\'erations 
$_{(n)},\ n\geq 0$, sur $L{\cal A}$ de fa\c{c}on \'evidente. 

Ensuite, on d\'efinit $U{\cal A}$ comme le quotient de l'enveloppe vertex 
de $L{\cal A}$ (cf. {\em op. cit.}, 10 (a)) par l'id\'eal vertex engendr\'e 
par les \'el\'ements $a\opmun x - ax,\ a\in A, x\in A\cup V_2$.    

Donc, on aura  
$$
(U{\cal{A}})_0 = A,\ (U{\cal{A}})_1 = 0,\ 
(U{\cal{A}})_2 = V_2,\ 
$$
$$
(U{\cal{A}})_3 = \partial V_2,\ 
(U{\cal{A}})_4 = \partial^2V_2\oplus S^2V_2,
$$     
\begin{center}
  etc.
\end{center}
  
\vspace{0.25in}

On veut maintenant voir qu'il est possible de considérer 
les algèbres de Frobenius comme des cas particuliers des alg\'ebro\"\i des 
de Virasoro. 

%\begin{definition}
%Une $\mathbb{C}$-algèbre $F$ sera dite de Frobenius si:
%\begin{enumerate}
%\item c'est une algèbre associative, commutative et unitaire (d'unité $e$);
%\item elle est munie d'une forme bilinéaire, symétrique et non-dégénérée
%$$\fonctionsansnom{F\times F}{\mathbb{C}}{(x,y)}{\langle x,y\rangle}$$
%\end{enumerate} qui est invariante dans le sens suivant:
%$$\langle x.y,z\rangle=\langle x,y.z\rangle.$$
%\end{definition}

Soit $F$ une algèbre de Frobenius sur $k$. On définit 
alors les trois opérations $\opun$, $\opzt$ et $\optr$ comme suit:
\begin{eqnarray*}
x\opun y &:=& 2xy,\\
x\opzt y &:=& xy,\\
\text{et }x\optr y &:=& \langle x,y\rangle.
\end{eqnarray*}
Il reste alors à s'assurer que ces opérations vérifient bien toutes les propriétés requises par la definition d'une algébroïde de Virasoro. Mais compte tenu de l'associativité et de la commutativité du 
produit dans la définition de l'algèbre de Frobenius, toutes ces 
vérifications sont immédiates. Donc $F$ devient une alg\'ebro\"\i de de 
Virasoro sur $k$; on la d\'esigne par $Vir(F)$. L'alg\`ebre vertex enveloppante 
$UVir(F)$ sera not\'ee $V_F$. 

\vspace{0.25in}

Explicitement, $V_F$ est d\'etermin\'ee par la propri\'et\'e universelle
suivante: pour une alg\`ebre vertex $V$ sur $k$ quelconque, un morphisme 
d'alg\`ebres vertex $V_F \lra V$ est la m\^eme chose qu'un morphisme 
de $k$-modules $\phi:\ F \lra V$ satisfaisant les conditions suivantes: 
$$  
\phi(x)\opun\phi(y) = 2\phi(xy)
$$
$$
\phi(x)\opz\phi(y) = \partial\phi(xy)
$$
$$
\phi(x)\optr\phi(y) = \langle x,y\rangle\cdot 1_V
$$

\paragraph{Rôle de l'élément unité} 
 
On va montrer maintenant que l'unité de l'algèbre de Frobenius initiale  
joue un rôle particulier dans notre nouvelle structure: en effet, 
elle peut être vue comme un vecteur de Virasoro dans $V_F$. 

On commence par rappeler la définition d'un tel vecteur: 

\begin{definition}
  Un vecteur de Virasoro de charge $c\in k$ dans une algèbre vertex  
  $\BZ_{\geq 0}$-gradu\'ee $V$ est 
  un élément $L\ \in\ V_2$ tel que 
\begin{gather}
L_{(0)}=\partial \label{cond1}\\
L_{(1)|V_j}=jId;\label{cond2}\\
L_{(2)}L=0;\label{cond3}\\
L_{(3)}L=\frac{c}{2}1.
\end{gather}
\end{definition}

\begin{proposition}
  L'élément unité $e$ de l'algèbre de Frobenius $F$ 
  est un vecteur de Virasoro dans $V = V_F$.
\end{proposition}

\begin{proof}
{\bfseries \'Equation \ref{cond1}}
Si $a$ est un élément de $V_0$, alors $e\opz a\in V_1$ donc comme $V_1$ est nul, $e\opz a=0$ et par conséquent, $e\opz$ coïncide avec $\partial$ sur $V_0$.

Si $x\in V_2$, alors
\begin{eqnarray*}
  e_{(0)}x & = & \partial \partial^{-1}(e_{(0)}x)\\
  & = & \partial (e\opzt x)\\
  & = & \partial (e.x)\\
  & = & \partial x.
\end{eqnarray*}
Dans ce cas, l'application $e\opz$ coïncide encore avec $\partial$.

Enfin, si $x\in V_3$, alors on peut trouver un élément $y$ dans $V_2$ tel que $x=\partial y$; on a alors
\begin{eqnarray*}
  e\opz x & = & e\opz \partial y\\
  & = & \partial(e\opz y) - \partial e\opz y\\
  & = & \partial(\partial y)+0.e\opmun y\\
  & = & \partial x.
\end{eqnarray*}
Encore une fois, l'égalité recherchée est vérifiée.

On a donc bien $e\opz=\partial$.

{\bfseries \'Equation \ref{cond2}}
Soit $a$ un élément de $V_0$; on sait que $e\opun a=a\opun e$ et une des remarques initiales nous permet alors d'affirmer que $e\opun a=0$. La condition recherchée est donc bien vérifiée pour $V_0$.

Soit $x\ \in\ V_2$. Par construction, on a $e\opun x=2e.x=2x$ et donc $e_{(1)|V_2}=2Id$.

Soit $x\ \in\ V_3$. Il existe $y\ \in\ V_2$ tel que $\partial y=x$. D'où:
\begin{eqnarray*}
  e\opun x & = & e\opun\partial y\\
  & = & \partial(e\opun y) - \partial e\opun y\\
  & = & \partial(2y)+e\opz y\\
  & = & 2\partial y + \partial y\\
  & = & 3\partial y.
\end{eqnarray*}

Au final, on a bien $e_{(1)|V_j}=jId$.

{\bfseries \'Equation \ref{cond3}}
Cette condition est trivialement vérifiée car $e_{(2)}e$ est 
un élément de $V_1$. 
\end{proof}

\vspace{0.5in} 

\centerline{\bf Bibliographie}

\vspace{0.25in}

[GMS] V.Gorbounov, F.Malikov, V.Schechtman, Gerbes of chiral differential 
operators. II. Vertex algebroids, {\em Inv. Math.}, {\bf 155} 
(2004), 605 - 680. 

[K] V.Kac, Vertex algebras for the beginners, University Lecture Series, 
{\bf 10}, American Mathematical Society, Providence, Rhode Island, 1997.

\end{document}